\documentclass[a4paper,11pt]{article}
\usepackage{amsfonts,amssymb,amsmath, dsfont,relsize,accents}
\usepackage{bm}
\usepackage{mathrsfs}  
\usepackage{color}
\usepackage{upgreek}
\usepackage[english]{babel}
\usepackage{graphicx}
\usepackage{microtype}
\usepackage[colorlinks=true, pdfstartview=FitV, linkcolor=blue, citecolor=blue, urlcolor=blue,pagebackref=false]{hyperref}
\usepackage{graphicx}
\usepackage{epstopdf}

\definecolor{labelkey}{gray}{.8}
\definecolor{refkey}{gray}{.8}

\definecolor{darkred}{rgb}{0.9,0.1,0.1}

 \makeatletter
 \@addtoreset{equation}{section}
 \makeatother

 \makeatletter
 \@addtoreset{enunciato}{section}
 \makeatother

 \newcounter{enunciato}[section]

 \newtheorem{ittheorem}{Theorem}
 \newtheorem{itlemma}{Lemma}
 \newtheorem{itproposition}{Proposition}
 \newtheorem{itcorollary}{Corollary}
 \newtheorem{itdefinition}{Definition}
 \newtheorem{itremark}{Remark}
 \newtheorem{itclaim}{Claim}
 \newtheorem{itfact}{Fact}
 \newtheorem{itconjecture}{Conjecture}

 \newenvironment{theorem}{\addtocounter{enunciato}{1}
 \begin{ittheorem}}{\end{ittheorem}}

 \newenvironment{lemma}{\addtocounter{enunciato}{1}
 \begin{itlemma}}{\end{itlemma}}

 \newenvironment{proposition}{\addtocounter{enunciato}{1}
 \begin{itproposition}}{\end{itproposition}}

 \newenvironment{corollary}{\addtocounter{enunciato}{1}
 \begin{itcorollary}}{\end{itcorollary}}

 \newenvironment{definition}{\addtocounter{enunciato}{1}
 \begin{itdefinition}}{\end{itdefinition}}

 \newenvironment{remark}{\addtocounter{enunciato}{1}
 \begin{itremark}}{\end{itremark}}

 \newenvironment{claim}{\addtocounter{enunciato}{1}
 \begin{itclaim}}{\end{itclaim}}

 \newenvironment{fact}{\addtocounter{enunciato}{1}
 \begin{itfact}}{\end{itfact}}

 \newenvironment{conjecture}{\addtocounter{enunciato}{1}
 \begin{itconjecture}}{\end{itconjecture}}

 \newcommand{\be}[1]{\begin{equation}\label{#1}}
 \newcommand{\ee}{\end{equation}}

 \newcommand{\bl}[1]{\begin{lemma}\label{#1}}
 \newcommand{\el}{\end{lemma}}

 \newcommand{\br}[1]{\begin{remark}\label{#1}}
 \newcommand{\er}{\end{remark}}

 \newcommand{\bt}[1]{\begin{theorem}\label{#1}}
 \newcommand{\et}{\end{theorem}}

 \newcommand{\bd}[1]{\begin{definition}\label{#1}}
 \newcommand{\ed}{\end{definition}}

 \newcommand{\bcl}[1]{\begin{claim}\label{#1}}
 \newcommand{\ecl}{\end{claim}}

 \newcommand{\bfact}[1]{\begin{fact}\label{#1}}
 \newcommand{\efact}{\end{fact}}

 \newcommand{\bp}[1]{\begin{proposition}\label{#1}}
 \newcommand{\ep}{\end{proposition}}

 \newcommand{\bc}[1]{\begin{corollary}\label{#1}}
 \newcommand{\ec}{\end{corollary}}

 \newcommand{\bcj}[1]{\begin{conjecture}\label{#1}}
 \newcommand{\ecj}{\end{conjecture}}

 \newcommand{\bpr}{\begin{proof}}
 \newcommand{\epr}{\end{proof}}

 \newcommand{\bprlem}[1]{\begin{proofof}{\it Lemma \ref{#1}}.\,\,}
 \newcommand{\eprlem}{\end{proofof}}

 \newcommand{\bprthm}[1]{\begin{proofof}{\it Theorem \ref{#1}}.\,\,}
 \newcommand{\eprthm}{\end{proofof}}

 \newcommand{\bprprop}[1]{\begin{proofof}{\it Proposition \ref{#1}}.\,\,}
 \newcommand{\eprprop}{\end{proofof}}

 \newcommand{\bi}{\begin{itemize}}
 \newcommand{\ei}{\end{itemize}}

 \newcommand{\ben}{\begin{enumerate}}
 \newcommand{\een}{\end{enumerate}}

 \newenvironment{proof}{\noindent {\em Proof}.\,\,}{\hspace*{\fill}$\halmos$\medskip}
 \newenvironment{proofof}{\noindent {\em Proof of\,\,}}{\hspace*{\fill}$\halmos$\medskip}
 \newcommand{\halmos}{\rule{1ex}{1.4ex}}

 \newcommand{\one}{{\mathchoice {1\mskip-4mu\mathrm l}
         {1\mskip-4mu\mathrm l}
         {1\mskip-4.5mu\mathrm l}
         {1\mskip-5mu\mathrm l}}}

\def \E {{\mathbb E}}

\def \N {{\mathbb N}}
\def \P {{\mathbb P}}

\def \R {{\mathbb R}}
\def \Z {{\mathbb Z}}
\def \T {{\mathbb T}}

\def \lra \leftrightarrow

\def \ra {\rightarrow}
\def \ba {\begin{array}}
\def \ea {\end{array}}

\def \lra {\longrightarrow}

\def \T {{\mathbb{T}}}
\def \lra {{\leftrightarrow}}

\def \ov {\overline}

\def \Ll {\left}
\def \Rr {\right}
\def \subset {\subseteq}

\def\one{\rlap{\mbox{\small\rm 1}}\kern.15em 1}

\newlength{\dhatheight}

\begin{document}
\title{Improved asymptotic estimates for the contact process with stirring}

\author{Anna Levit\textsuperscript{1}, Daniel Valesin\textsuperscript{2}}
\footnotetext[1]{Department of Mathematics, University of British Columbia. \url{anna.levit@math.ubc.ca}}
\footnotetext[2]{Johann Bernoulli Institute, University of Groningen. \url{d.rodrigues.valesin@rug.nl}}
\date{September 14, 2015}
\maketitle

\begin{abstract}
We study the contact process with stirring on $\mathbb{Z}^d$. In this process, particles occupy vertices of $\mathbb{Z}^d$; each particle dies with rate 1 and generates a new particle at a randomly chosen neighboring vertex with rate $\lambda$, provided the chosen vertex is empty. Additionally, particles move according to a symmetric exclusion process with rate $N$. For any $d$ and $N$, there exists $\lambda_c$ such that, when the system starts from a single particle, particles go extinct when $\lambda < \lambda_c$ and have a chance of being present for all times when $\lambda > \lambda_c$. Durrett and Neuhauser  proved that  $\lambda_c$ converges to 1 as $N$ goes to infinity, and Konno, Katori and Berezin and Mytnik obtained dimension-dependent asymptotics for this convergence, which are sharp in dimensions 3 and higher. We obtain a lower bound which 
is new in dimension 2 and also gives the sharp asymptotics in dimensions 3 and higher. Our proof involves an estimate for two-type renewal processes which is of independent interest.\\

\noindent \textit{Keywords:} interacting particle systems, contact process, contact process with rapid stirring
\end{abstract}

\section{Introduction}
The \textit{contact process on $\Z^d$ with birth rate $\lambda > 0$ and stirring rate $N > 0$} is the Markov process $(\xi_t)_{t \geq 0}$ on $\{0,1\}^{\Z^d}$ with generator $\mathcal{L} = \mathcal{L}^{(N)}_\text{exc}+ \mathcal{L}^{(\lambda)}_\text{cont}$ defined, for any function $f:\{0,1\}^{\Z^d} \to \R$ that only depends on finitely many coordinates and any $\xi \in \{0,1\}^{\Z^d}$, by:
\begin{align}
&(\mathcal{L}^{(N)}_\text{exc}f)(\xi) = N\cdot \sum_{\substack{\{x,y\} \subset \Z^d:\\x\sim y}} \left(f(\xi^{x \leftrightarrow y}) - f(\xi) \right), \text{ where }\xi^{x\leftrightarrow y}(z) = \begin{cases} \xi(y) &\text{if } z = x,\\ \xi(x) &\text{if } z = y,\\\xi(z) &\text{otherwise},\end{cases}\\
&\nonumber (\mathcal{L}^{(\lambda)}_\text{cont}f)(\xi) = \sum_{\substack{x \in \Z^d:\\\xi(x) = 1}} \left(f(\xi^{x \leftarrow 0}) - f(\xi) \right)  +\frac{\lambda}{2d}\cdot \sum_{\substack{x \in \Z^d:\\\xi(x) = 1}}\sum_{\substack{y \in \Z^d:\\y \sim x}} \left(f(\xi^{y \leftarrow 1}) - f(\xi)\right),\\&\hspace{6cm}\text{where }\xi^{x \leftarrow i}(z) = \begin{cases} i &\text{if } z = x,\\\xi(z)&\text{otherwise,}\end{cases}\quad i\in \{0,1\}
\end{align}
and we write $x \sim y$ if $|x - y|_1 = 1$. We interpret each vertex $x\in\Z^d$ as a region of space which can either be empty (represented by state 0) or contain a particle (state 1). The dynamics can be described as follows:
\begin{itemize}
\item $\mathcal{L}_\text{exc}$ describes the motion of particles, which is according to an \textit{exclusion dynamics} with rate $N$. This means that for each pair $x \sim y$ with $\xi(x) = 1$ and $\xi(y) = 0$, the particle at $x$ jumps to $y$ with rate $N$. For each pair $x \sim y$ with $\xi(x) = \xi(y)= 1$, it will be useful to think that the particles at $x$ and $y$ exchange positions with rate $N$ (though this amounts to no change in the configuration $\xi$).
\item $\mathcal{L}_\text{cont}$ describes the birth and death of particles, which is according to a \textit{contact process} with rate $\lambda$, described as follows. Each particle dies with rate 1 and gives birth with rate $\lambda$; when a birth occurs, a new particle is placed at a position chosen uniformly at random among the neighbors of the parent (if the chosen position is not empty, no new particle is born).
\end{itemize}

This process is then a mixture of two well-studied classes of interacting particle systems: the exclusion process, introduced by Spitzer in \cite{spitzer}, and the contact process, introduced by Harris in \cite{harris}. It was first considered in \cite{durrneu} by Durrett in Neuhauser, who were motivated by the earlier work \cite{mfp} by DeMasi, Ferrari and Lebowitz. In \cite{mfp}, the authors added exclusion dynamics to particles on Glauber-type spin systems and showed that, if the rate of the exclusion is taken to infinity, the system converges to the solution of an associated reaction-diffusion equation.

Our interest will be to consider the contact process with stirring from the point of view of its extinction-survival phase transition. To explain what we mean by this, let us first consider the basic contact process (that is, take the above definition with $N = 0$). Assume $\xi_0 = \mathds{1}_{\{0\}}$ (the indicator function of the origin) and consider the probability (which depends on $d$ and $\lambda$)
\begin{equation}\P\Ll[\text{for all $t$ there exists } x\text{ such that }\xi_t(x) = 1\Rr].\label{eq:phase_transition}\end{equation}
This is non-decreasing in $\lambda$. We say that the process \textit{dies out} if \eqref{eq:phase_transition} is zero and \textit{survives} if it is positive. The phase transition for the contact process is the statement that there exists $\lambda_c(\Z^d) \in (0,\infty)$ such that the process survives if and only if $\lambda > \lambda_c(\Z^d)$ (for a proof of this, and of all the facts we state here about the basic contact process, see \cite{lig99}).

Since attempted births on occupied vertices (which we call ``collisions'') produce no new particles, it is easy to show that the process given by the total number of particles, $(\sum_{x \in \Z^d} \xi_t(x))_{t \geq 0}$, is stochastically dominated by a branching process with birth rate $\lambda$ and death rate $1$. This comparison yields $\lambda_c(\Z^d) \geq 1$, and in fact it is known that
$$\frac{2d}{2d-1} \leq \lambda_c(\Z^d) \leq 2,\qquad d \geq 1.$$

When stirring is introduced, one can consider $\lambda_c(\Z^d, N)$, still defined as the infimum of the values of $\lambda$ for which \eqref{eq:phase_transition} is positive. In \cite{durrneu}, Durrett and Neuhauser proved that
\begin{equation}\label{eq:phase_transition_stirring}\lim_{N \to \infty} \lambda_c(\Z^d, N) = 1,\qquad d \geq 1.\end{equation}
This result roughly means that, as we take $N \to \infty$ -- thus allowing particles to move more and more between each birth or death event -- collisions have less effect and the critical rate approaches that of the associated branching process. 

Regarding the rate of convergence in \eqref{eq:phase_transition_stirring}, Konno showed in \cite{konno} that for every $d$ there exist positive constants $c_d$, $C_d$ such that, for all $N \geq 1$,
\begin{equation}\label{eq:thm_konno}
\lambda_c(\Z^d, N) - 1 \in \begin{cases} \left[\frac{c_d}{N},\;\frac{C_d}{N}\right]&\text{if } d\geq 3;\\[.2cm]\left[\frac{c_d \log N}{N},\; \frac{C_d \log N}{N}\right]&\text{if } d = 2;\\[.2cm] \left[\frac{c_d}{N^{1/3}},\;\frac{C_d}{N^{1/3}}\right]&\text{if } d = 1.\end{cases}
\end{equation}
For $d \geq 3$, in \cite{katori}, Katori improved this to
\begin{equation}\label{eq:thm_katori}
\frac{1}{2d(2d-1)} < \liminf_{N\to\infty} N(\lambda_c(\Z^d, N) - 1)< \limsup_{N\to\infty} N(\lambda_c(\Z^d, N) - 1) < \frac{G(0,0)-1}{2d},
\end{equation}
where $G(\cdot, \cdot)$ is the Green function of discrete-time, simple random walk on $\Z^d$; that is, if $\beta_1, \beta_2,\ldots$ are independent and uniformly distributed on $\{z: z \sim 0\}$, then $$G(x,y) = \mathds{1}_{\{x=y\}} +\E [\#\{n \geq 1: \Sigma_{i=1}^n \beta_i = y -x \}],$$
where $\#$ denotes the cardinality of a set. More recently, Berezin and Mytnik proved that, for $d \geq 3$, the upper bound in \eqref{eq:thm_katori} is sharp:
\begin{equation} \label{eq:thm_bm}
\lim_{N \to \infty} N\cdot (\lambda_c(\Z^d,N) - 1) = \frac{G(0,0)-1}{2d}.\end{equation}

We now state our main results. For a measurable set $A \subset \R$, $|A|$ denotes the Lebesgue measure of $A$.
\begin{theorem} \label{thm:main} Assume $d \geq 2$ and let $(A_s, B_s)_{s \geq 0}$ denote the positions of two particles under exclusion dynamics on $\Z^d$ with rate 1 (per edge) and $A_0 \sim B_0$. 
Suppose $(\lambda_N)_{N \geq 1}$ is a sequence satisfying
$$\limsup_{N\to\infty} \frac{\lambda_N - 1}{\frac{1}{dN}\cdot \E\Ll[|\{t \leq N: A_t \sim B_t\}|\Rr]} < 1.$$
Then, if $N$ is large enough, the contact process with birth rate $\lambda_N$ and stirring rate $N$ dies out.
\end{theorem}
The following helps relate this theorem to the bounds given earlier:
\begin{proposition}\label{prop:exc_rw} Let $A_s, B_s$ be as in Theorem \ref{thm:main}. Let $A_s'$ and $B_s'$ be independent, continuous-time simple random walks on $\Z^d$ which jump from each vertex with rate $2d$, and $A_0' \sim B_0'$.  
If $d \geq 3$, then 
\begin{equation}\label{eq:local_d3}
\E\Ll[|\{t < \infty: A_t \sim B_t\}|\Rr] =\E\Ll[|\{t < \infty: A_t' \sim B_t'\}|\Rr] =  \frac{G(0,0)-1}{2}.
\end{equation}
If $d = 2$, then
\begin{equation}\label{eq:local_d2}\lim_{t \to \infty} \frac{\E\Ll[|\{s < t: A_t \sim B_t\}|\Rr]}{\log t} =\lim_{t \to \infty} \frac{\E\Ll[|\{s < t: A'_t \sim B'_t\}|\Rr]}{\log t} = \frac{1}{2\pi}.\end{equation}
\end{proposition}
\eqref{eq:local_d3} already appeared in \cite{romanleonid}. \eqref{eq:local_d2}, which presents a higher technical challenge, is new. By putting together Proposition \ref{prop:exc_rw} and Theorem \ref{thm:main}, we re-obtain the sharp bound of Berezin and Mytnik for $d \geq 3$ and obtain a new bound for $d = 2$:
\begin{corollary}
If $d \geq 3$, then 
$$\liminf_{N \to \infty} N\cdot (\lambda_c(\Z^d,N) - 1) \geq \frac{G(0,0)-1}{2d}.$$
If $d = 2$, then
\begin{equation}\label{eq:result_d2}\liminf_{N \to \infty} \frac{N}{\log N} \cdot (\lambda_c(\Z^d,N) - 1) \geq \frac{1}{4\pi}.\end{equation}
\end{corollary}
Our bound for $d=2$ is of the correct order given in \eqref{eq:thm_konno}. No matching upper bound has been proved, but since both cases $d =  2$ and $d \geq 3$ are treated by the unified statement of Theorem \ref{thm:main}, and since the resulting bound is sharp in the second case, it is reasonable to expect that \eqref{eq:result_d2} is also sharp.

In our proof of Proposition  \ref{prop:exc_rw}, we study the processes $(A_s - B_s)_{s \geq 0}$ and $(A'_s - B'_s)_{s \geq 0}$. The latter is exactly a random walk on $\Z^d$ and the former is the same random walk except that its jump rates at the set of neighbors of the origin are modified, and in particular it can never reach the origin. Hence, denoting by $\mathcal{N}_0$ the set of neighbors of 0,
$$(\mathds{1}{\{A_s - B_s \in \mathcal{N}_0\}})_{s \geq 0} \text{ and } (\mathds{1}{\{A_s' - B_s' \in \mathcal{N}_0\cup \{0\}\}})_{s \geq 0}, $$
are two-type (0 and 1) renewal processes, and the amounts of time they spend at stage 0 have the same distribution. 

Motivated by this, we consider the following more general setting. Let $U^{(1)}$, $U^{(2)}$ and $V$ be positive random variables and consider three independent sequences of independent random variables,
\begin{equation}\label{eq2:def_uvs}
U^{(1)}_0, U^{(1)}_1,U^{(1)}_2\ldots \sim U^{(1)},\qquad U^{(2)}_0, U^{(2)}_1,U^{(2)}_2\ldots \sim U^{(2)},\qquad V_0,V_1,V_2,\ldots \sim V. 
\end{equation}
For $i \in \{1,2\}$, define sequences $(S^{(i)}_n)_{n \geq 0}$ as follows:
\begin{equation}\label{eq:def_ss}
S^{(i)}_0 \equiv 0\text{ and, for all } n \geq 0,\; S^{(i)}_{2n+1}-S^{(i)}_{2n} = U^{(i)}_n \text{ and } S^{(i)}_{2n+2}-S^{(i)}_{2n+1} = V_n.
\end{equation}
Also define, for $i \in \{1,2\}$ and $t > 0$,
\begin{equation}\label{def_kappa}
\kappa_t^{(i)} =  \left|[0,t] \cap \left(\bigcup_{n=0}^\infty [S_{2n}^{(i)},\;S_{2n+1}^{(i)}]\right)\right|= t - \left|[0,t] \cap \left(\bigcup_{n=0}^\infty [S_{2n+1}^{(i)}, S_{2n+2}^{(i)}] \right) \right|.
\end{equation}
We think of the sequence $(S^{(i)}_n)_{n \geq 0}$ as describing a renewal process which alternates between a ``$u^{(i)}$-state'' (where it stays for an amount of time distributed as $U^{(i)}$) and a ``$v$-state'' (where it stays for an amount of time distributed as $V$). Then, $\kappa_t^{(i)}$ represents the total amount of time spent in the $u^{(i)}$-state before instant $t$. We prove:
\begin{theorem}\label{thm:renewal}
If $\mathbb{E}[(U^{(1)})^2]$, $\mathbb{E}[(U^{(2)})^2] < \infty$ and $\mathbb{E}[V] = \infty$,
then
\begin{equation}
\lim_{t\to \infty} \frac{\E[\kappa^{(1)}_t]}{\E[\kappa^{(2)}_t]} = \frac{\E[U^{(1)}]}{\E[U^{(2)}]}.
\end{equation}
\end{theorem}
We think this result could be useful in other settings, particularly in other estimates involving local times of exclusion processes with finitely many particles.

The paper is organized as follows. In Section \ref{s:construct}, we give a construction of the contact process with stirring that allows us to separately consider the genealogy of particles and their motion. The construction is also very convenient to compare the process, particle by particle, with the associated branching process that bounds it from above. Although this construction and comparison were already described with words (and somewhat vaguely) in \cite{durrneu} and \cite{romanleonid}, as far as we know this is the first time that they are given explicitly. In Section \ref{s:proofs}, we show how Theorem \ref{thm:renewal} implies Proposition \ref{prop:exc_rw} and then use the construction of Section \ref{s:construct} to prove Theorem \ref{thm:main}. In Section \ref{s:renewal}, we prove Theorem \ref{thm:renewal}.

\section{Construction of coupled processes}
\label{s:construct}

\textbf{Basic genealogical process.} We start by giving a construction of a continuous-time branching process which will be useful for coupling, and then comparing, processes of interest. Our construction will depend on the parameter $\lambda > 0$.

Let $\T$ be the tree defined as follows. The vertex set of $\T$ (which, by abuse of notation, is also denoted $\T$), is $\{o\} \cup \left(\cup_{n=1}^\infty \mathbb{N}^n\right)$, where $o$ is a distinguished element called the \textit{root}. The edge set is
$$E(\mathbb{T}) = \{\{o,i\}:i\in\N \}\cup\{\{(i_1,\ldots,i_n),(i_1,\ldots,i_n,i_{n+1})\}:i_1,\ldots,i_{n+1}\in\N\}.$$
In case $\alpha = (i_1,\ldots,i_n)$ and $\beta = (i_1,\ldots, i_n, i_{n+1})$, we say that $\alpha$ is the \textit{parent} of $\beta$ (and denote this by $\alpha = p(\beta)$) and $\beta$ is the $i_{n+1}$-th \textit{child} of $\alpha$ (denoted $\beta = s_{i_{n+1}}(\alpha)$). The same terminology and notation is used in case $\alpha = o$ and $\beta = i \in \N$.

Let $H_\T = \{B_\alpha,\;D_\alpha\}_{\alpha \in \T}$ be a family of independent Poisson point processes on $[0,\infty)$ so that each $B_\alpha$ has rate $\lambda$ and each $D_\alpha$ has rate $1$. We regard each $B_\alpha$ or $D_\alpha$ as a random discrete subset of $[0,\infty)$.

We will define a process $(\Psi_t)_{t \geq 0}$ with state space $\{0,1,-1\}^\T$. Our terminology will be as follows: each vertex in $\T$ is called a \textit{particle}; a particle is \textit{present} if it is in state 1 and \textit{absent} if in state 0 or $-1$. We start setting $\Psi_0 = \mathds{1}_{\{o\}}$. Now assume $\Psi_t$ has already been defined up to time $t \geq 0$. Let
$$t' = \inf\{s > t: s \in B_\alpha \cup D_\alpha \text{ for some $\alpha$ with } \Psi_t(\alpha) = 1\}.$$
We set $\Psi_s = \Psi_t$ for every $s \in (t, t')$; in particular, this means that, if $\Psi_t(\alpha) \neq 1$ for every $\alpha$, then $t'= \infty$ and $\Psi_s = \Psi_t$ for every $s > t$. In case $t'<\infty$, define $\Psi_{t'}$ as follows. For the unique $\alpha$ such that $t'\in B_\alpha \cup D_\alpha$, consider the two cases:
\begin{itemize}
\item if $t'\in D_\alpha$, we set $\Psi_{t'}(\alpha) = -1$ and $\Psi_{t'}(\beta) = \Psi_t(\beta)$ for all $\beta \neq \alpha$;
\item if $t' \in B_\alpha$, let $n$ be the smallest natural number for which $\Psi_t(s_n(\alpha))=0$, and set $\Psi_{t'}(s_n(\alpha)) = 1$ and $\Psi_{t'}(\beta) = \Psi_t(\beta)$ for all $\beta \neq s_n(\alpha)$. 
\end{itemize}
It should be clear that, with this construction, the number of present particles at time $t$, $\#\{\alpha:\Psi_t(\alpha) = 1\}$, is a continuous-time branching process with birth rate equal to $\lambda$ and death rate equal to 1. This consideration also shows that the above prescription defines the process $(\Psi_t)$ for all $t \geq 0$, that is, no finite-time explosion occurs in the application of the recursive procedure.\\[.2cm]
\textbf{Particle positions.} In addition to $\lambda$, we now also fix $N > 0$. 

We now take two more random objects, independent of each other and independent of $H_\T$:
\begin{itemize}
\item a family $H_{\mathbb{Z}^d} = \{L_e: e\in E(\mathbb{Z}^d)\}$ ($E(\mathbb{Z}^d)$ denotes the set of nearest neighbors edges of $\mathbb{Z}^d$). Each $L_e$ is a Poisson point processes with rate $1$ on $[0,\infty)$, and these processes are all independent;
\item a family $M = \{M_\alpha(t): \alpha \in \T\backslash\{o\},\; t \geq 0\}$ of independent random vectors of $\mathbb{Z}^d$ (that is, if $(\alpha_1,t_1),\ldots, (\alpha_k, t_k)$ are all distinct, then $M_{\alpha_1}(t_1),\ldots, M_{\alpha_k}(t_k)$ are independent). Each $M_\alpha(t)$ is uniformly distributed on the set of the neighbors of the origin of $\mathbb{Z}^d$.
\end{itemize}
Given a realization of $H_{\mathbb{Z}^d}$, we define a function $\rho: \{(x,s,t): x \in \mathbb{Z}^d,\; 0\leq s \leq t\} \to \mathbb{Z}^d$, which we call a \textit{flow}, as follows: for each $x \in \mathbb{Z}^d$ and $s \geq 0$, $t \mapsto \rho(x,s,t)$ is the unique function that is constant by parts and satisfies $\rho(x,s,s) = x$ and, for all $t \geq s$,
\begin{equation}\begin{array}{ll}\rho(x,s,t) = \rho(x,s,t-) \text{ if } t \notin \cup_{z: z\sim \rho(x,s,t-)}\;L_{\{\rho(x,s,t-),z\}};\\[.2cm]
\rho(x,s,t) = z \neq \rho(x,s,t-)  \text{ if } t \in L_{\{\rho(x,s,t-),z\}}\end{array}
\end{equation}
(we hence view $t \in L_{\{x,y\}}$ as an ``order'' either to jump from $x$ to $y$ or to jump from $y$ to $x$, and $\rho$ is the path obtained by following all the orders that are encountered). We also let
$$\rho_N(x,s,t) = \rho(x,Ns,Nt).$$

For each $\alpha \in \T$, let
$$\uptau^-_\alpha = \inf\{t: \Psi_t(\alpha) = 1\},\qquad \uptau^+_\alpha = \sup\{t:\Psi_t(\alpha) = 1\}.$$
Then define
$$X_o(t) = \begin{cases}\rho_N(0,0,t)&\text{if } 0\leq t < \uptau^+_o;\\ \triangle&\text{otherwise},\end{cases}$$
where $\triangle$ denotes a ``cemetery'' state. Now assume that processes $(X_\alpha(t))_{t \geq 0}$ on $\mathbb{Z}^d \cup \{\triangle\}$ have been defined for all $\alpha \in \bigcup_{m=0}^n \;\N^m$, and that these processes satisfy $X_\alpha(t) \neq \triangle$ if and only if $\Psi_t(\alpha) =1$. Fix $\alpha \in \N^{n+1}$. In case $\uptau^-_\alpha = \infty$, put $X_\alpha(t) = \triangle$ for all $t$. Otherwise,
\begin{equation} \label{eq:appear}X_\alpha(t) = \begin{cases}\rho_N\left(X_{p(\alpha)}(\uptau^-_\alpha) + M_\alpha(\uptau^-_\alpha), \uptau^-_\alpha, t\right),&\text{if } t \in [\uptau^-_\alpha, \uptau^+_\alpha);\\[.2cm]\triangle&\text{otherwise.} \end{cases} \end{equation}
This inductively defines $X_\alpha(t)$ for all $\alpha \in \T$ and $t \geq 0$, so that $X_\alpha(t) \neq \triangle$ if and only if $\Psi_t(\alpha) = 1$. We call $X_\alpha(t)$ the \textit{location} of particle $\alpha$ at time $t$ (with the understanding that absent particles are located in the cemetery state). The definition \eqref{eq:appear} thus means that when a particle appears, it is placed on a location obtained as a random neighbor of its parent's location at the time, and then it moves according to the flow $\rho_N$ until it disappears.
\begin{remark}
Define the process
\begin{equation*}  \psi_t(x) = \#\{\alpha \in \T: \Psi_t(\alpha) = 1,\; X_\alpha(t) = x\}, \qquad x \in \mathbb{Z}^d,\; t \geq 0. \end{equation*}
Although we will not need it in the sequel, it is instructive to discuss its behavior at this point. In this process, particles occupy positions in $\mathbb{Z}^d$; it is possible that any number of particles occupy a single position. Each particle disappears with rate 1 and gives birth at a randomly chosen neighboring position with rate $\lambda$ (births are not forbidden at occupied sites). Moreover, edges of $\Z^d$ contain jump instructions which are Poisson($N$) clocks; the effect of a jump instruction at $\{x,y\}$ is that \textbf{all} particles from $x$ jump to $y$ simultaneously, and vice-versa. This is the process studied by Katori in \cite{katori}, inspired by the ``binary contact path process'' of Griffeath (\cite{griffeath}). Our task now is to construct, in this same probability space, the contact process with stirring $(\xi_t)_{t\geq 0}$: this will amount, quite simply, to keeping track of particle positions and forbidding births at occupied sites. With this construction, at any point in time, the set of particles present in $\xi_t$ will be a subset of the set of particles in $\psi_t$.
\end{remark}

\noindent \textbf{Contact process with stirring.} We will now construct the contact process with stirring on $\mathbb{Z}^d$ (denoted $(\xi_t)_{t \geq 0}$, with state space $\{0,1\}^{\mathbb{Z}^d}$) and its underlying genealogical process ($(\Xi_t)_{t \geq 0}$ on $\{0,1,-1\}^\T$). The construction of these processes will depend on $(H_\T, H_{\mathbb{Z}^d}, M)$ that have been used above and on the paths $\{(X_\alpha(t))_{t \geq 0}: \alpha \in \T\}$ that have been defined (hence  $\lambda$ and $N$ are fixed throughout).

We will give a recursive definition of $\Xi$ that mimics that of $\Psi$ and will set
\begin{equation} \label{eq:def_xi} \xi_t(x) = \#\{\alpha \in \T: \Xi_t(\alpha) = 1,\; X_\alpha(t) = x\}, \qquad x \in \mathbb{Z}^d,\; t \geq 0. \end{equation}
The construction will guarantee that, for all $s \geq 0$,
\begin{align}\label{eq:good_xi_range} &\xi_s(x) \in \{0,1\} \text{ for all }s \geq 0,\;x \in \mathbb{Z}^d;\\[.2cm]
\label{eq:coup_comp}&\{\alpha: \Xi_s(\alpha) = 1\} \subseteq \{\alpha:\Psi_s(\alpha) = 1\}; \\[.2cm]
&\label{eq:fact2}
\text{if }\Xi_s(\alpha) = 1, \text{ then } \inf\{n:\Xi_s(s_n(\alpha)) = 0\} = \inf\{n: \Psi_s(s_n(\alpha)) = 0\}.
\end{align}

We set $\Xi_0 = \mathds{1}_{\{o\}}$. Assume that $(\Xi_s)_{0 \leq s \leq t}$ has already been defined up to time $t \geq 0$ (and hence, by \eqref{eq:def_xi}, $(\xi_s)_{0 \leq s \leq t}$ is also defined) and that they satisfy \eqref{eq:coup_comp}, \eqref{eq:fact2} and \eqref{eq:good_xi_range} for all $s \leq t$.
We let
$$t' = \inf\{s > t: s \in B_\alpha \cup D_\alpha \text{ for some $\alpha$ with } \Xi_t(\alpha) = 1\}.$$
For $s \in (t, t')$, we let $\Xi_s = \Xi_t$ (note that this, together with \eqref{eq:def_xi}, implies that $\xi$ is then defined in $[0, t')$). If $t' < \infty$, for the unique $\alpha$ for which $t'\in B_\alpha \cup D_\alpha$, we consider the cases:
\begin{itemize}
\item if $t'\in D_\alpha$, set $\Xi_{t'}(\alpha) = -1$ and $\Xi_{t'}(\beta) = \Xi_t(\beta)$ for all $\beta \neq \alpha$;
\item if $t' \in B_\alpha$, let $n$ be the smallest natural number for which $\Xi_t(s_n(\alpha))=0$, and let $\beta = s_n(\alpha)$. We consider two further cases:
\begin{itemize}
\item if $\xi_{t'-}(X_\alpha(t'-) + M_\beta(t'-)) = 0$ (that is, if the position where $\beta$ is supposed to appear is empty), set $\Xi_{t'}(\beta) = 1$ and $\Xi_{t'}(\gamma) = \Xi_{t}(\gamma)$ for all $\gamma \neq \beta$;
\item if $\xi_{t'-}(X_\alpha(t'-) + M_\beta(t'-)) \neq 0$ (that is, if that position is occupied), set $\Xi_{t'}(\beta) = -1$ and $\Xi_{t'}(\gamma) = \Xi_{t}(\gamma)$ for all $\gamma \neq \beta$.
\end{itemize}
\end{itemize}
These rules define $\Xi$ and $\xi$ for all times. Each step in the induction preserves \eqref{eq:coup_comp}, \eqref{eq:fact2} and \eqref{eq:good_xi_range}, so they are satisfied for all times. Also note that
\begin{equation}
\begin{array}{ll}\text{if } \Xi_s(\alpha) = 1 \text{ for some } s,\text{ then } &\inf\{t: \Xi_t(\alpha) = 1\} = \inf\{t: \Psi_t(\alpha) = 1\} =  \uptau^-_\alpha \text{ and }\\[.2cm] &\sup\{t:\Xi_t(\alpha) = 1\} = \sup\{t:\Psi_t(\alpha) = 1\}= \uptau^+_\alpha,\end{array}
\end{equation}
that is, if a particle $\alpha$ ever becomes present in $\Xi$, then the set of times for which it is present is the same for $\Xi$ and $\Psi$. Also,
\begin{equation}\begin{array}{l}
\text{if, for } \alpha, \beta \in \T,\; \uptau_\alpha^- < \uptau_\beta^- < \infty \text{ and } X_\alpha(\uptau^-_\beta-) = X_{p(\beta)}(\uptau_\beta^--)+ M_\beta(\uptau_\beta^--),\\[.2cm]\hspace{6cm} \text{ then } \Xi_t(\beta) \neq 1 \text{ for all } t\geq 0, 
\end{array}\end{equation}
that is, a particle $\beta$ cannot appear at position $x \in \Z^d$ and time $t$ if $x$ is occupied at $t$ by a particle $\alpha$ that appeared earlier than $t$.

For $t \geq 0$, we denote by $\mathfrak{F}_t$ the smallest sigma-algebra under which the point processes $H_\T \cap [0,t],\; H_{\mathbb{Z}^d} \cap [0,t]$ and the processes  $\{(M_\alpha(s))_{0\leq s \leq t}:\alpha \in \T\}$ are measurable. Note that $(\Psi_t, \Xi_t, \xi_t)_{t\geq 0}$ is then adapted to $(\mathfrak{F}_t)_{t \geq 0}$.

\section{Proof of main result}
\label{s:proofs}
Given $\alpha \in \T$ and $t \geq 0$ so that $\Xi_t(\alpha) = 1$, let $n$ be the smallest integer so that $\Xi_t(s_n(\alpha)) = 0$ (or equivalently, by \eqref{eq:fact2}, so that $\Psi_t(s_n(\alpha))= 0$). Then consider the subgraph of $\mathbb{T}$ given by 
$$\hat \T(\alpha, t) = \{\alpha\} \cup \{\beta: \text{the unique path in $\T$ from $\beta$ to $o$ intersects $\{s_m(\alpha): m \geq n\}$}\}$$
(together with all edges that are incident to two vertices in the above set). We call $\hat \T(\alpha,t)$ the set of \textit{fresh descendants of $\alpha$ at time $t$}. It is a connected subtree of $\T$. We also define, for $\alpha$ and $t$ such that $\Xi_t(\alpha) = 1$ and $t'> t$,
$$\mathcal{N}^\Xi(\alpha,t,t') = \#\{\beta \in \hat \T(\alpha,t):\Xi_{t'}(\beta) = 1\},\quad \mathcal{N}^\Psi(\alpha,t,t') = \#\{\beta \in \hat \T(\alpha,t):\Psi_{t'}(\beta) = 1\}.$$
We observe that, for all $\alpha, t$ with $\Xi_t(\alpha) = 1$,
\begin{equation}\label{eq:NleN}\mathcal{N}^\Xi(\alpha,t,t') \leq  \mathcal{N}^\Psi(\alpha,t,t')\end{equation}
and, since $\left(\mathcal{N}^\Psi(\alpha,t,t+s) \right)_{s \geq 0}$ is a branching process with birth rate $\lambda$ and death rate 1, we have
\begin{equation}
\label{eq:Psi_bran}
\text{on } \{\Xi_t(\alpha) = 1 \},\text{ for all } t'\geq t,\; \mathbb{E}\left[\mathcal{N}^\Psi(\alpha,t,t') \mid\mathfrak{F}_t\right] = e^{(\lambda - 1)(t'-t)}.
\end{equation}

We now define two classes of events, $I(\alpha,t)$ and $J(\alpha,t)$, which will both correspond to situations where there is a $0\to 1$ transition in $\Psi$ which does not occur in $\Xi$. We will be able to guarantee that, if one of these events occurs, the number of  particles present in $\Xi$ is strictly smaller than that in $\Psi$. Both definitions will depend on the following quantity:
\begin{equation} t^* = \frac{1}{\log N}. \label{eq:def_t_star}\end{equation}

We start with the definition of $I(\alpha,t)$. Fix $\alpha \in \T$, $t \geq 0$ and let $n$ be the smallest integer such that $\Xi_t(s_n(\alpha)) = 0$. Let $\beta = s_n(\alpha)$ and $\gamma = s_{n+1}(\alpha)$. Assume that the following occurs:
\begin{equation}\label{eq:event_i1}
\Xi_t(\alpha) = 1,\qquad
(D_\alpha \cup D_\beta \cup D_\gamma \cup B_\beta \cup B_\gamma) \cap [t,t+t^*] = \varnothing,\qquad
\#(B_\alpha \cap [t,t+t^*]) = 2.
\end{equation}
Additionally, letting $B_\alpha \cap [t,t+t^*] = \{T_1,T_2\}$ with $T_1 < T_2$, assume that
\begin{equation}\label{eq:event_i2}
X_\alpha(T_2-) + M_\gamma(T_2-) = X_\beta(T_2-).
\end{equation}
(This implies that either $\beta$ obstructs the appearance of $\gamma$, or some earlier particle obstructs the appearance of both $\beta$ and $\gamma$). Let $I(\alpha, t)$ be the event described in \eqref{eq:event_i1} and \eqref{eq:event_i2}. We then have
\begin{equation}\label{eq:cons_i}
\text{on } I(\alpha,t),\;\mathcal{N}^\Xi(\alpha,t,t+ t^*) \leq 2 \text{ and }\mathcal{N}^\Psi(\alpha,t,t+t^*) = 3.\end{equation}

We now turn to the similar definition of $J(\alpha,t)$. Again fix $\alpha \in \T$, $t \geq 0$, $n$ the smallest integer such that $\Xi_t(s_n(\alpha)) = 0$ and $\beta = s_n(\alpha)$. Also let $\gamma' = s_1(\beta)$. Assume that
\begin{equation}\label{eq:event_j1}\begin{split}
\Xi_t(\alpha) = 1,\qquad
&(D_\alpha \cup D_\beta  \cup D_{\gamma'} \cup B_{\gamma'}) \cap [t,t+t^*] = \varnothing,\\
&\#(B_\alpha \cap [t,t+t^*]) = \#(B_\beta \cap [t,t+t^*]) = 1.\end{split}
\end{equation}
Letting $B_\alpha \cap [t, t+t^*] = \{T_1\}$ and $B_\beta \cap [t, t+t^*] = \{T_2\}$, we will also require that
\begin{equation}\label{eq:event_j2} 
T_1 < T_2,\qquad X_\alpha(T_2-) = X_\beta(T_2-) + M_{\gamma'}(T_2-).
\end{equation}
(This implies that either $\beta$ obstructs the appearance of $\gamma'$, or some earlier particle obstructs the appearance of both $\beta$ and $\gamma'$). We define $J(\alpha,t)$ as the event specified by all the requirements in \eqref{eq:event_j1} and \eqref{eq:event_j2}, and note that
\begin{equation}
\label{eq:cons_j}
\text{on } J(\alpha,t),\; \mathcal{N}^\Xi(\alpha,t,t+t^*) \leq 2 \text{ and } \mathcal{N}^\Psi(\alpha,t,t+t^*) = 3.
\end{equation}

Due to the invariance under time shift of $H_\T, H_{\mathbb{Z}^d}$ and $M$,
\begin{equation}
\label{eq:time_shift} \text{on } \{\Xi_t(\alpha) = 1\},\quad\;\P\left[I(\alpha,t)\;|\;\mathfrak{F}_t\right] = \P[I(o,0)]\quad \text{ and }\quad \P\left[J(\alpha,t)\;|\;\mathfrak{F}_t\right] = \P[J(o,0)].
\end{equation}
Also, inspecting the above definitions one can show that
\begin{equation}
\label{eq:I_equal_J}\P\left[I(o,0) \right] = \P\left[J(o,0)\right].
\end{equation}

We are now ready for the key estimate of this section. For any $t \geq 0$,
\begin{align*}
&\E\left[\left. \#\{\alpha \in \T: \Xi_{t+t^*}(\alpha) = 1\}\;\right|\;\mathfrak{F}_t\right] =\E\left[\left.\sum_{\substack{\alpha \in \T: \Xi_t(\alpha) = 1}}\mathcal{N}^\Xi(\alpha, t, t+t^*) \right|\;\mathfrak{F}_t\right]\\[.2cm]
&\stackrel{\eqref{eq:NleN},\eqref{eq:cons_i},\eqref{eq:cons_j}}{\leq}\E\left[\left.\sum_{\substack{\alpha \in \T: \Xi_t(\alpha) = 1}}\mathcal{N}^\Psi(\alpha, t, t+t^*) \right|\;\mathfrak{F}_t\right] - \E\left[\left. \sum_{\alpha \in \T: \Xi_t(\alpha) = 1} \left(\mathds{1}_{I(\alpha,t)} + \mathds{1}_{J(\alpha,t)}\right)\right|\;\mathfrak{F}_t\right]\\[.2cm]
&\hspace{.5cm}\stackrel{\eqref{eq:Psi_bran},\eqref{eq:time_shift}}{=} \#\{\alpha \in \T:\Xi_t(\alpha) = 1\} \cdot (e^{t^*(\lambda - 1)} - \P[I(o,0)] - \P[J(o,0)])\\[.2cm]
&\hspace{.8cm} \stackrel{\eqref{eq:I_equal_J}}{=} \#\{\alpha \in \T:\Xi_t(\alpha) = 1\} \cdot (e^{t^*(\lambda - 1)} - 2\P[I(o,0)]).
\end{align*}
Applying this estimate recursively, for any $k \in \N$ we have
$$\begin{aligned}\mathbb{E}\left[\#\{x \in \mathbb{Z}^d: \xi_{kt^*}(x) = 1\}\right] &\stackrel{\eqref{eq:def_xi}}{=} \E\left[\#\{\alpha \in \T:\Xi_{kt^*}(\alpha) = 1\}\right] \leq (e^{t^*(\lambda - 1)} - 2\P[I(o,0)] )^k.\end{aligned}$$
This shows that
\begin{equation}\label{eq:criterion}
\begin{split} &\text{if } e^{t^*(\lambda - 1)} - 2\P[I(o,0)]  < 1, \text{ then }\\[.2cm]&\qquad \P[\xi_t \neq \varnothing \; \text{ for all } t] = \lim_{k \to \infty} \P[\xi_{kt^*} \neq \varnothing] \leq \lim_{k \to \infty} \E[\#\{x: \xi_{kt^*}(x) = 1\}]= 0.
\end{split}
\end{equation}

We now estimate $\P[I(o,0)]$. We denote $\beta = s_1(o)$, $\gamma = s_2(o)$ and 
$$E = \{(D_o \cup D_\beta\cup D_\gamma\cup B_\beta \cup B_\gamma) \cap [0,t^*] = \varnothing,\; \#(B_o \cap [0,t^*]) = 2\},$$
so that
\begin{equation}\label{eq:probE}
\P[E] = \frac{(\lambda t^*)^2}{2} \exp\{-3t^* - 3\lambda t^*\}.\end{equation}
On $E$, let $T_1 < T_2$ denote the two elements of $B_o \cap [0,t^*]$. Let $f$ denote the joint density function of $T_1, T_2$ conditioned on $E$; then, $f(t_1,t_2) = \frac{2}{(t^*)^2}\cdot \mathds{1}_{\{0 < t_1< t_2 < t^*\}}.$ Now,
\begin{align*}
&\P[I(o,0)\mid E] =  \sum_{\substack{z_1,z_2 \in \mathbb{Z}^d:\\z_1,z_2 \sim 0}}\P[M_\beta(T_1) = z_1,\;M_\gamma(T_2) = z_2]  \\ 
&\times \sum_{x,y\in\Z^d}\int_0^{t^*} \int_{t_1}^{t^*} \frac{2}{(t^*)^2}\;\P\left[\rho_N(0,0,t_1) = x,\;\rho_N(0,0,t_2) = y,\; \rho_N(x + z_1,t_1,t_2) = y + z_2 \right]\mathrm{d}t_2\mathrm{d}t_1 \\[.2cm]
&\hspace{1.0cm}=\frac{1}{2d} \cdot \sum_{\substack{z_1 \in \mathbb{Z}^d:\\z_1 \sim 0}}\P[M_\beta(T_1) = z_1]\cdot \int_0^{t^*} \int_{t_1}^{t^*} \frac{2}{(t^*)^2}\cdot \P[\rho_N(0,0,t_2 - t_1)  \sim \rho_N(z_1,0,t_2-t_1)] \; \mathrm{d}t_2\mathrm{d}t_1\\[.2cm]
&\hspace{1.0cm}=\frac{1}{2d} \cdot  \int_0^{t^*} \int_{t_1}^{t^*} \frac{2}{(t^*)^2}\cdot \P[\rho_N(0,0,t_2 - t_1)  \sim \rho_N(\bar z,0,t_2-t_1)] \; \mathrm{d}t_2\mathrm{d}t_1
\end{align*}
for any fixed $\bar z \sim 0$ in $\mathbb{Z}^d$. Changing variables twice, this is further equal to
\begin{align}
\nonumber&\frac{1}{2d}  \cdot \frac{2}{(t^*)^2}\int_0^{t^*} \P\left[\rho_N(0,0,s) \sim \rho_N(\bar{z},0,s)\right]\cdot (t^* - s) \;\mathrm{d}s\\[.2cm]
&=\frac{1}{dt^*N}\int_0^{t^*N}\P\left[\rho_1(0,0,u) \sim \rho_1(\bar{z},0,u)\right]\cdot \left(1 - \frac{u}{t^*N}\right) \;\mathrm{d}u.\label{eq:after_change}
\end{align}
We now fix $\varepsilon \in (0,1)$ whose value will be chosen later. The expression in \eqref{eq:after_change} is larger than
\begin{align}
\nonumber&\frac{1}{dt^*N} \cdot \left(1 - \varepsilon\right)\int_0^{\varepsilon t^*N}\P\left[\rho_1(0,0,u) \sim \rho_1(\bar{z},0,u)\right] \;\mathrm{d}u.\\[.2cm]
&=\frac{1-\varepsilon}{dt^*N} \cdot \E\left[|\{u \leq \varepsilon t^*N: \rho_1(0,0,u) \sim \rho_1(\bar{z},0,u)\}|\right].\label{eq:estimate_exp_as}\end{align}
We abbreviate
$$g(t) = \E[|\{u \leq t: \rho_1(0,0,u) \sim \rho_1(\bar{z},0,u)\}|];$$
by Proposition \ref{prop:exc_rw}, to be proved in the next subsection, we have $\lim_{t\to\infty} \frac{g(\varepsilon t^*N)}{g(N)} = 1$ for any $d \geq 2$. Putting together this fact, \eqref{eq:local_d3}, \eqref{eq:local_d2}, \eqref{eq:probE} and  \eqref{eq:estimate_exp_as} we see that, if $N$ is large enough (depending on $\varepsilon$),
\begin{equation}
\nonumber\P[I(o,0)] \geq  \frac{(1-\varepsilon)^2}{dt^*N} \cdot g(N) \cdot \frac{(\lambda t^*)^2}{2} \cdot \exp\{-3t^* -3\lambda t^* \}.
\end{equation}
If we further assume that $\lambda = \lambda_N \xrightarrow{N \to \infty} 1$ and use \eqref{eq:def_t_star}, then for $N$ large enough we have
\begin{equation}
\label{eq:final_est_I}
\P[I(o,0)] \geq (1-\varepsilon)^3 \cdot \frac{t^*  g(N)}{2dN}.
\end{equation}

We are now ready to conclude. Assume that $\lambda = \lambda_N = 1 + \theta \frac{g(N)}{dN}$ for some $\theta < 1$ (so in particular, by by Proposition \ref{prop:exc_rw}, we have $\lambda_N \xrightarrow{N \to \infty} 1$). Then, if $N$ is large enough,
$$e^{t^*(\lambda - 1)} - 1 \leq (1+\varepsilon)\cdot \theta \frac{t^*g(N)}{dN} \stackrel{(\star)}{<} (1-\varepsilon)^3 \cdot  \frac{t^*g(N)}{dN} \stackrel{\eqref{eq:final_est_I}}{<} 2 \P[I(o,0)],$$
where the inequality ($\star$) holds if $\varepsilon$ is small enough that $\frac{(1-\varepsilon)^3}{1+\varepsilon} > \theta$. The desired result now follows from \eqref{eq:criterion}.



\subsection{Exclusion dynamics of two particles: proof of Proposition \ref{prop:exc_rw}}
We will denote the vectors in the canonical basis of $\Z^d$ by $\vec{e}_1,\ldots,\vec{e}_d$. We write
$$\mathcal{N}_0 = \{ \vec{e}_1,-\vec{e}_1,\ldots, \vec{e}_d, -\vec{e}_d \},\quad \ov{\mathcal{N}}_0 = \mathcal{N}_0 \cup \{0\}.$$
For $d \geq 3$, as in the Introduction, we denote by $G(x,y)$ the Green function of discrete-time, simple random walk on $\Z^d$ and by $(A_s,B_s)_{s \geq 0}$ the process given by the positions of two particles moving on $\Z^d$ under exclusion dynamics with rate 1 per edge of $\Z^d$. This means that $(A_s,B_s)_{s \geq 0}$ is the Markov process on $(\Z^d)^2$ with generator
$$\begin{aligned}(\mathcal{L}_\text{exc}f)(a,b) =& \sum_{z \in \mathcal{N}_0}  \left[ \mathds{1}_{\{a+z \neq b\}}\cdot  \left( f(a+z,b) - f(a,b)\right) +\mathds{1}_{\{b+z \neq a\}}\cdot  \left(f(a,b+z) - f(a,b)\right)\right. \\&\hspace{6cm}\left.+ \mathds{1}_{\{a + z  = b\} \cup \{b+z = a\}} \cdot \left(f(b,a) - f(a,b)\right)\right].\end{aligned}$$
We always assume that $A_0 - B_0 \in \mathcal{N}_0$, that is, the particles are initially at neighboring positions. 

Now let
$$X_s = A_s - B_s,\;s \geq 0;$$
observe that $X_s \in \Z^d\backslash \{0\}$ for all $s$ and that $X_0 \in \mathcal{N}_0$. The generator of $(X_s)$ is 
$$(\mathcal{L}_Xf)(x) = \begin{cases}2\sum_{z \in \mathcal{N}_0} (f(x + z) - f(x)) &\text{if } x \notin \mathcal{N}_0;\\[.2cm](f(-x) - f(x))+  2\sum_{z \in \mathcal{N}_0\backslash \{-x\}} (f(x+z) - f(x))  &\text{otherwise.}  \end{cases}$$
For comparison purposes, it will be useful to also consider the process $(Y_s)_{s\geq 0}$, which has generator
$$(\mathcal{L}_Yf)(x) = 2\sum_{z \in \mathcal{N}_0} (f(x + z) - f(x));$$
we also assume $Y_0 \in \mathcal{N}_0$. $(Y_s)$ is hence a continuous-time simple random walk on $\Z^d$ that jumps from each vertex with total rate $4d$. $(X_s)$ and $(Y_s)$ have the same behavior except when $(X_s)$ is in some position $z\in\mathcal{N}_0$, in which case, instead of jumping to 0 with rate 2, it jumps to $-z$ with rate 1.

Now let
$$\begin{aligned}&U^X = \inf\{t \geq 0: X_t \notin \mathcal{N}_0\};\\
&U^Y = \inf\{t \geq 0: Y_t \notin \ov{\mathcal{N}}_0\};\\
&U^{Y,\mathcal{N}_0} = |\{t \in [0, U^Y]: Y_t \in \mathcal{N}_0\};\quad U^{Y,0} = |\{t \in [0, U^Y]: Y_t = 0 \},
\end{aligned}$$
so that $U^Y = U^{Y,\mathcal{N}_0} + U^{Y,0}$. By a simple computation involving geometrically distributed random variables, it is easy to check (or see the proof of Lemma 3.4 in \cite{romanleonid}) that
\begin{equation} \label{eq:form_leo}
\E\left[U^X\right] = \E\left[U^{Y,\mathcal{N}_0}\right] = \frac{1}{4d-2},\qquad \E\left[U^Y\right] = \frac{2d+1}{8d^2-4d}.
\end{equation}
\begin{proofof}{\textit{Proposition \ref{prop:exc_rw}}.} First assume $d \geq 3$. Let $q$ denote the probability that $(X_s)$ never returns to $\mathcal{N}_0$ after first leaving it. Note that this is the same as the probability that $(Y_s)$ never returns to $\ov{\mathcal{N}}_0$ after first leaving it. We then have
$$\begin{aligned}&\E\left[\left|\{s < \infty: X_s \in \mathcal{N}_0\}\right|\right] \\[.2cm]&= \frac{\mathbb{E}[U^X]}{q} \stackrel{\eqref{eq:form_leo}}{=}   \frac{\mathbb{E}[U^{Y,\mathcal{N}_0}]}{q} = \E\left[\left|\{s < \infty: Y_s \in \mathcal{N}_0\}\right|\right]=\frac{1}{4d}\sum_{z\in\mathcal{N}_0} G(0,z )= \frac12(G(0,0)-1),\end{aligned}$$
where the last equality follows from $G(0,0) = 1 + \frac{1}{2d}\sum_{z\in\mathcal{N}_0} G(z,0) = 1 + G(0,\vec{e}_1)$, by symmetry of the Green function.

For $d = 2$ we have
\begin{equation}\label{eq:final_lim}\lim_{t \to \infty} \frac{\E\left[|\{s\leq t: X_s \in \mathcal{N}_0\}|\right]}{\log t} = \lim_{t\to\infty} \frac{\E\left[|\{s\leq t: Y_s \in \ov{\mathcal{N}}_0\}|\right]}{\log t}\cdot \frac{\E\left[|\{s\leq t: X_s \in \mathcal{N}_0\}|\right]}{\E\left[|\{s\leq t: Y_s \in \ov{\mathcal{N}}_0\}|\right]}.\end{equation}
Now, Theorem \ref{thm:renewal} implies that
\begin{equation}\label{eq:final_lim1}
\lim_{t\to\infty} \frac{\E\left[|\{s\leq t: X_s \in \mathcal{N}_0\}|\right]}{\E\left[|\{s\leq t: Y_s \in \ov{\mathcal{N}}_0\}|\right]} =  \frac{\E\left[U^X\right]}{\E\left[U^Y\right]}\stackrel{\eqref{eq:form_leo}}{=} \frac{2d}{2d+1} = \frac{4}{5}.
\end{equation} 
Now let us treat the first quotient in \eqref{eq:final_lim}. By the Local Central Limit Theorem (Theorem 2.5.6 in \cite{lawler}), for any $z \in \Z^d$,
$$\frac{\P[Y_s = z]}{f_{4s}(z)} \xrightarrow{s \to \infty} 1,$$
where $f_{t}(x) = \frac{1}{2\pi t}e^{-|x|^2/(2t)}$ is the probability density function of a Gaussian vector $(V,W)$ so that $V$ and $W$ are independent and have mean 0 and variance $t$. Hence,
\begin{equation}\label{eq:final_lim2}\frac{\E\left[|\{s\leq t: Y_s \in \ov{\mathcal{N}}_0\}|\right]}{\log t} = \frac{1}{\log t}\sum_{z \in \ov{\mathcal{N}}_0} \int_0^t \P\left[Y_s = z\right]\;\mathrm{d}s \xrightarrow{t \to \infty} \frac{2d+1}{2\pi\cdot 4} = \frac{5}{8\pi}\end{equation}
and similarly,
\begin{equation}
\label{eq:final_lim3} \frac{\E\left[|\{s\leq t: Y_s \in \mathcal{N}_0\}|\right]}{\log t}\xrightarrow{t \to \infty} \frac{2d}{2\pi \cdot 4} = \frac{1}{2\pi}.
\end{equation}
Now, \eqref{eq:final_lim}, \eqref{eq:final_lim1}, \eqref{eq:final_lim2} and \eqref{eq:final_lim3}  imply the desired result.
\end{proofof}

\section{Two-type renewal processes: proof of Theorem \ref{thm:renewal}}
\label{s:renewal}
Take $(V_n)_{n\geq 0}$, $(U^{(i)}_n)_{n \geq 0}$, $(S^{(i)}_n)_{n\geq 0}$ and $\kappa^{(i)}_t$, for $i = 1,2$, as in \eqref{eq2:def_uvs}, \eqref{eq:def_ss}, \eqref{def_kappa}. Theorem \ref{thm:renewal} follows immediately from the two lemmas:
\begin{lemma}\label{lem:det_non_det}
If $U^{(1)} \equiv 1$, $\E[U^{(2)}] = 1$  and $\E[(U^{(2)})^2] < \infty$, then
\begin{equation}
\lim_{t\to \infty} \frac{\E[\kappa^{(1)}_t]}{\E[\kappa^{(2)}_t]} = 1.
\end{equation}
\end{lemma}
\begin{lemma}\label{lem:det_det}
Let $\alpha^{(1)},\; \alpha^{(2)} > 0$. If $U^{(1)} \equiv \alpha^{(1)}$, $U^{(2)} \equiv \alpha^{(2)}$ and $\E[V] = \infty$, then
\begin{equation}\label{eq:des_det_det}
\lim_{t\to \infty} \frac{\E[\kappa^{(1)}_t]}{\E[\kappa^{(2)}_t]} = \frac{\alpha^{(1)}}{\alpha^{(2)}}.
\end{equation}
\end{lemma}
Before we prove these results, let us give a definition. For $i \in \{1,2\}$ and $t \geq 0$, we let $N^{(i)}_t$ be the unique value of $n$ such that $S^{(i)}_{2n} \leq t < S^{(i)}_{2n + 2}$; note that
\begin{align}
\label{eq:nt} \kappa^{(i)}_t &= \sum_{n=0}^{N^{(i)}_t - 1} U^{(i)}_n + \min\left(t,\; S^{(i)}_{2N^{(i)}_t + 1}\right) - S^{(i)}_{2N^{(i)}_t}\\[.2cm]
&= t - \sum_{n=0}^{N^{(i)}_t -1} V_n - \max\left(0,\;t - S^{(i)}_{2N^{(i)}_t + 1} \right).\label{eq:nt2}
\end{align}
(Here and in what follows, we interpret sums of the form $\sum_{n=a}^b$ with $a>b$ to be equal to zero). We also observe that, since $U^{(i)}$ and $V$ are assumed to be positive random variables, we have
\begin{equation}\label{eq:n_goes_inf}
\lim_{\substack{t \to \infty\\\text{a.s.}}}N^{(i)}_t = \infty\quad  \text{ and } \quad\lim_{t \to \infty} \E[N^{(i)}_t] = \infty.
\end{equation}
We also claim that
\begin{equation}\label{eq:n_t_goes_zero}
\text{if } \E[U^{(i)} + V] = \infty, \text{ then } \lim_{t\to\infty} \frac{1}{t} \E[N^{(i)}_t] = 0.
\end{equation}
Indeed, the fact that $N_t/t \to 0$ almost surely as $t \to \infty$ follows from the Law of Large Numbers, so it is sufficient to prove that $\sup_{t \geq 1} \E[(N_t/t)^2] < \infty$. This can be shown by fixing $\varepsilon >0$, $\delta > 0$ such that $\P[U^{(i)} + V > \delta] > \varepsilon$ and noting that
$$N_t^{(i)} \leq \min\left\{m \geq 0:\sum_{n=0}^{m-1} \mathds{1}\{U^{(i)}_n + V_n > \delta\} > \lceil t/\delta \rceil\right\};$$
the right-hand side is distributed as a sum of $\lceil t/\delta \rceil$ geometric random variables with parameter larger than $\varepsilon$; this gives $\mathbb{E}[(N_t^{(i)})^2] \leq \left\lceil \frac{t}{\delta}\right\rceil\cdot \frac{2-\varepsilon}{\varepsilon}$. Hence \eqref{eq:n_t_goes_zero} is proved.

\begin{proofof}{\textit{Lemma \ref{lem:det_non_det}}.} Define
$$\begin{aligned}\Delta_m = &\max\left\{|t - t'|: \sum_{n=0}^{m-1} U^{(1)}_n \leq t \leq \sum_{n=0}^{m} U^{(1)}_n,\; \sum_{n=0}^{m-1} U^{(2)}_n \leq t' \leq \sum_{n=0}^{m} U^{(2)}_n\right\}, \qquad m \geq 0.
\end{aligned}$$

We will now prove two claims.
\begin{claim} There exists $C > 0$ such that, for all $k \geq 1$,
\begin{equation} \label{eq:eqlem}\mathbb{E}\left[\max_{0\leq m \leq k} \Delta_m \right] \leq Ck^{\frac34}.\end{equation}
\end{claim}
To prove this, start noting that, since $U^{(1)}\equiv 1$, $\Delta_m$ is equal to
\begin{align*}
\max\left\{\left|m-\sum_{n=0}^{m}U^{(2)}_n\right|, \left|m+1-\sum_{n=0}^{m-1}U^{(2)}_n\right|\right\}\leq \max\left\{\left|m+1-\sum_{n=0}^{m}U^{(2)}_n\right|, \left|m-\sum_{n=0}^{m-1}U^{(2)}_n\right|\right\}+1,
\end{align*}
so that
\begin{equation*}
\max_{0\leq m\leq k}\Delta_m\leq  \max_{1\leq m\leq k+1}\left|m -\sum_{n=0}^{m-1}U^{(2)}_n\right|+1.
\end{equation*}
Since the $U^{(2)}_n$ are independent and identically distributed with mean 1, the Reflection Principle (\cite{lawler}, Proposition 1.6.2) implies that
\begin{equation*}
\mathbb{E}\left[\frac{\max_{1\leq m\leq k+1}\Delta_m}{k^{3/4}}\right] \leq 2 \cdot \mathbb{E} \left[\frac{\left|k+1-\sum_{n=0}^{k}U^{(2)}_n\right|}{k^{3/4}}\right]+ \frac{1}{k^{3/4}}.
\end{equation*}
Write $W_k = \left|k+1-\sum_{n=0}^{k}U^{(2)}_n\right|$. By the law of the iterated logarithm,
\begin{equation} \limsup_{\substack{k\to\infty\\\text{a.s.}}} \frac{W_k}{\sqrt{k \log \log k}} < \infty \qquad \text{hence} \qquad \label{anna}  \lim_{\substack{k\to\infty\\\text{a.s.}}}\;\frac{W_k}{k^{3/4}} = 0.\end{equation}
Then, 
\begin{align*}
\mathbb{E}\left[\frac{W_k}{k^{3/4}}\right] & \leq \mathbb{E}\left[\frac{W_k}{k^{3/4}} \cdot \mathds{1}_{\left\{\frac{W_k}{k^{3/4}} \leq 1\right\}}\right] +  \mathbb{E}\left[\left(\frac{W_k}{k^{3/4}} \right)^2\right]\\[.2cm]
&=\mathbb{E}\left[\frac{W_k}{k^{3/4}} \cdot \mathds{1}_{\left\{\frac{W_k}{k^{3/4}} \leq 1\right\}}\right] +  \frac{\text{Var}(W_k)}{k^{3/2}}=\mathbb{E}\left[\frac{W_k}{k^{3/4}} \cdot \mathds{1}_{\left\{\frac{W_k}{k^{3/4}} \leq 1\right\}}\right] +  \frac{k+1}{k^{3/2}}\cdot \text{Var}(U^{(2)}).
\end{align*}
The expectation on the right-hand side vanishes as $k \to \infty$ by \eqref{anna} and the dominated convergence theorem. This completes the proof of the claim.

\begin{claim}
\begin{equation}\label{eq:bound}
\left| \kappa^{(2)}_t - \kappa^{(1)}_t\right| \leq \max_{0 \leq m \leq N^{(1)}_t} \Delta_m,\qquad t > 0.
\end{equation}
\end{claim}
To show this, fix $t > 0$. First consider the case $N^{(1)}_t = N^{(2)}_t$; then, it readily follows from \eqref{eq:nt} that $\left| \kappa^{(2)}_t - \kappa^{(1)}_t\right| \leq \Delta_{N^{(1)}_t - 1}$. Next, assume that 
\begin{equation} \label{eq:assume} N^{(1)}_t > N^{(2)}_t. \end{equation} We then have
$$\begin{aligned}\kappa^{(1)}_t &\stackrel{\eqref{eq:nt2}}{=}\;  t - \sum_{n=0}^{N^{(1)}_t -1} V_n - \max\left(0,\;t - S^{(1)}_{2N^{(1)}_t + 1} \right) \stackrel{\eqref{eq:assume}}{\leq}  t - \sum_{n=0}^{N^{(2)}_t } V_n  \stackrel{\eqref{eq:nt2}}{\leq} \kappa^{(2)}_t. \end{aligned}$$
Hence,
$$\begin{aligned}&\left| \kappa^{(2)}_t - \kappa^{(1)}_t\right| = \kappa^{(2)}_t - \kappa^{(1)}_t \stackrel{\eqref{eq:nt}}{=}\;\; \sum_{n=0}^{N^{(2)}_t - 1} U^{(2)}_n + \min\left(t,\; S^{(2)}_{2N^{(2)}_t + 1}\right) - S^{(2)}_{2N^{(2)}_t}\\
&\hspace{5.5cm}- \sum_{n=0}^{N^{(1)}_t - 1} U^{(1)}_n - \left(\min\left(t,\; S^{(1)}_{2N^{(1)}_t + 1}\right) - S^{(1)}_{2N^{(1)}_t}\right)\\[.2cm]
&\stackrel{\eqref{eq:assume}}{\leq}  \sum_{n=0}^{N^{(2)}_t - 1} U^{(2)}_n + \min\left(t,\; S^{(2)}_{2N^{(2)}_t + 1}\right) - S^{(2)}_{2N^{(2)}_t} - \sum_{n=0}^{N^{(2)}_t - 1} U^{(1)}_n\\[.2cm]
&\leq \Delta_{N^{(2)}_t - 1} \stackrel{\eqref{eq:assume}}{\leq} \max_{0 \leq m \leq N^{(1)}_t} \Delta_m. \end{aligned}$$
Now, a symmetric argument shows that, when $N^{(1)}_t < N^{(2)}_t$, we have
$$\left| \kappa^{(2)}_t - \kappa^{(1)}_t\right| \leq  \Delta_{N^{(1)}_t - 1} \stackrel{\eqref{eq:assume}}{\leq} \max_{0 \leq m \leq N^{(1)}_t} \Delta_m.$$
This completes the proof of \eqref{eq:bound}.

We are now ready to conclude. We have 
\begin{equation}\left|\frac{\mathbb{E}[\kappa^{(2)}_t]}{\mathbb{E}[\kappa^{(1)}_t]} - 1 \right|  = \mathbb{E}[\kappa^{(1)}_t]^{-1} \cdot \left|\mathbb{E}[\kappa^{(2)}_t] - \mathbb{E}[\kappa^{(1)}_t]\right| \leq \mathbb{E}[\kappa^{(1)}_t]^{-1}  \cdot \mathbb{E} \left[   \max_{0 \leq m \leq N^{(1)}_t} \Delta_m\right].\label{eq:partial}\end{equation}
Now note that $N^{(1)}_t$ only depends on $(V_n)_{n\geq 0}$ and $(\Delta_n)_{n \geq 0}$ only depends on  $(U^{(2)}_n)_{n\geq 0}$. Hence, $N^{(1)}_t$ is independent of $(\Delta_n)_{n\geq 0}$. Thus,
$$\begin{aligned}
&\mathbb{E} \left[   \max_{0 \leq m \leq N^{(1)}_t} \Delta_m\right] = \sum_{k=0}^\infty \mathbb{E}\left[\max_{0 \leq m \leq k} \Delta_m \;|\; N^{(1)}_t = k\right]\cdot \mathbb{P}\left[N^{(1)}_t = k\right]\\&\qquad
= \sum_{k=0}^\infty \mathbb{E}\left[\max_{0 \leq m \leq k} \Delta_m \right]\cdot \mathbb{P}\left[N^{(1)}_t = k\right] \stackrel{\eqref{eq:eqlem}}{\leq} C  \sum_{k=0}^\infty k^\frac34 \cdot  \mathbb{P}\left[N^{(1)}_t = k\right] = C\mathbb{E}\left[(N^{(1)}_t)^\frac34\right].
\end{aligned}$$
Since $N^{(1)}_t \to \infty$ almost surely as $t \to \infty$, we have $\mathbb{E}\left[(N^{(1)}_t)^\frac34\right]/\mathbb{E}\left[N^{(1)}_t\right] \to 0$
 as $t \to \infty$. This shows that the right-hand side of \eqref{eq:partial} vanishes as $t \to \infty$, completing the proof.
\end{proofof}

\begin{proofof}{\textit{Lemma \ref{lem:det_det}}.}
If $\alpha^{(1)} = \alpha^{(2)}$, the result follows from Lemma \ref{lem:det_non_det} and a change of scale. So we assume (without loss of generality) that $\alpha^{(1)} > \alpha^{(2)}$. It follows from the assumption that $U^{(i)} \equiv \alpha^{(i)}$ and \eqref{eq:nt} that
$$|\kappa^{(i)}_t - \alpha^{(i)} \cdot N^{(i)}_t | \leq \alpha^{(i)}, \qquad i \in \{1,2\},\; t \geq 0,$$
so, writing
$$\frac{\E[\kappa^{(1)}_t]}{\E[\kappa^{(2)}_t]} = \frac{\alpha^{(1)} \E[N^{(1)}_t] + \E[\kappa^{(1)}_t - \alpha^{(1)} N^{(1)}_t]}{\alpha^{(2)} \E[N^{(2)}_t] + \E[\kappa^{(2)}_t - \alpha^{(2)} N^{(2)}_t]}$$
and using \eqref{eq:n_goes_inf}, we see that \eqref{eq:des_det_det} will follow once we prove that
\begin{equation}
\label{eq:reduction_N}\lim_{t \to \infty} \frac{\E[N^{(1)}_t]}{\E[N^{(2)}_t]} = 1.
\end{equation}

With this in mind, we now proceed to bound $|N_t^{(2)} - N_t^{(1)}|$. Observe that
\begin{equation}
\label{eq:cond_n_maior}
N^{(i)}_t < k \qquad \text{ if and only if }\qquad  t < S^{(i)}_{2k} = \alpha^{(i)} k + \sum_{n=0}^{k-1} V_n.
\end{equation}
Since 
$$\alpha^{(2)} N^{(1)}_t + \sum_{n=0}^{N^{(1)}_t-1} V_n < \alpha^{(1)} N^{(1)}_t + \sum_{n=0}^{N^{(1)}_t-1} V_n =S_{2N^{(1)}_t} \leq t,$$
applying \eqref{eq:cond_n_maior} with $i = 2$ and $k = N^{(1)}_t$ implies that $N_t^{(2)} \geq N_t^{(1)}$.
We now claim that
\begin{equation}
\label{eq:claim_star_red}
N^{(2)}_t - N^{(1)}_t < 1+\inf \left\{ m \geq 0 : \alpha^{(2)} \cdot m + \sum_{n= N^{(1)}_t + 1}^{N^{(1)}_t + m } V_n > (\alpha^{(1)} - \alpha^{(2)}) (N^{(1)}_t + 1)\right\}.
\end{equation}
Indeed, assume $m$ belongs to the set of which the infimum is taken on the right-hand side. We have
\begin{align*}\alpha^{(2)} (N^{(1)}_t + m +1) + \sum_{n=0}^{N^{(1)}_t + m} V_n &= \alpha^{(2)} N^{(1)}_t + \sum_{n=0}^{N_t^{(1)}}V_n + \alpha^{(2)} + \left(\alpha^{(2)}m + \sum_{n= N^{(1)}_t +1}^{N^{(1)}_t + m } V_n\right)\\[.2cm]
&> \alpha^{(2)} N^{(1)}_t + \sum_{n=0}^{N^{(1)}_t} V_n + \alpha^{(2)} +  (\alpha^{(1)} - \alpha^{(2)}) (N^{(1)}_t+1)\\[.2cm]
& = \alpha^{(1)}( N^{(1)}_t+1) + \sum_{n=0}^{N^{(1)}_t} V_n  = S^{(1)}_{2N^{(1)}_t + 2} > t,
\end{align*}
so, applying \eqref{eq:cond_n_maior} with $i = 2$ and $k = N^{(1)}_t + m + 1$ we get $N_t^{(2)} < N_t^{(1)} + m + 1$, proving  \eqref{eq:claim_star_red}.

Now, if $V_0', V_1',\ldots$ are random variables distributed as $V$ and independent of $(U^{(i)}_n)_{n\geq 0}$ and $(V_n)_{n\geq 0}$, the right-hand side of \eqref{eq:claim_star_red} has the same distribution as
$$1+ \inf\left\{m \geq 0: \alpha^{(2)} m + \sum_{n=0}^{m-1} V_n' > (\alpha^{(1)} - \alpha^{(2)})( N^{(1)}_t + 1)\right\}.$$
Conditioning on $N_t^{(1)}$ (and denoting by $\mu_{N^{(1)}_t}$ the distribution of $N^{(1)}_t$), \eqref{eq:claim_star_red} then yields
\begin{align*}\E[|N^{(2)}_t - N^{(1)}_t|] &\leq \int_0^\infty  \E\left[1+\inf\left\{m \geq 0: \alpha^{(2)} m + \sum_{n=0}^{m-1} V_n' > (\alpha^{(1)} - \alpha^{(2)})(s+1) \right\}\right] \mu_{N^{(1)}_t} (ds)\\
&= \int_0^\infty \E\left[1+N^{(2)}_{(\alpha^{(1)} - \alpha^{(2)})(s+1)}\right] \; \mu_{N^{(1)}_t} (ds).
\end{align*}
Now fix $\varepsilon > 0$. By \eqref{eq:n_t_goes_zero}, there exists $C > 0$ such that $\E[N^{(2)}_s] \leq C + \varepsilon s \text{ for all } s \geq 0$, so the above gives
$$\E[|N^{(2)}_t - N^{(1)}_t|] \leq 1 + C + \varepsilon\cdot  \E[(\alpha^{(1)} - \alpha^{(2)})(N^{(1)}_t + 1)].$$
Finally, by \eqref{eq:n_goes_inf},
$$0 \leq \limsup_{t\to\infty} \frac{\E[|N^{(2)}_t - N^{(1)}_t|]}{\E[N^{(1)}_t]} \leq \varepsilon(\alpha^{(1)} - \alpha^{(2)}); $$
since $\varepsilon$ is arbitrary, this completes the proof of \eqref{eq:reduction_N}.
\end{proofof}

\begin{center}\textbf{Acknowledgments}\end{center}
\vspace{-.3cm}We would like to thank Leonid Mytnik for introducing us to the model and problems that we consider in this paper. We would also like to thank Ed Perkins for helpful discussions.

\end{document}